\documentclass[dvips, 10pt, a4paper]{article}
\usepackage{latexsym, amsmath}
\usepackage{amssymb}

\begin{document}

\begin{center}
{\LARGE Some inequalities on hemi-slant product submanifolds in a cosymplectic manifold}\\

\bigskip
\noindent
{\large Khushwant Singh, S. S. Bhatia }\\
\end{center}
\begin{quotation}
\noindent {\bf Abstract.} Recently, M. Atcken studied Contact
CR-warped product submanifolds in cosymplectic space forms and
established general sharp inequalities for CR-warped products in a
cosymplectic manifold [1]. In the present paper, we obtain an
inequality for the squared norm of the second fundamental form in
terms of constant $\phi-${\it sectional curvature} for hemi-slant
products in cosymplectic manifolds. An inequality for hemi-slant
warped products in a cosymplectc manifold is also given. The
equality case is considered.
\end{quotation}

\noindent
{\bf{M.S.C. 2000:}} 53C15, 53C40, 53C25

\section{Introduction}

Bishop and O' Neill [4] introduced the concept of warped products in
1969. They defined as follows\\

\noindent{\bf {Definition 1.1.}} Let $(B,g_B)$ and $(F,g_F)$ be two
Riemannian manifolds with Riemannian metric $g_B$ and $g_F$
respectively and $f$ a positive differentiable function on $B$. The
warped product $B\times_f F$ of $B$ and $F$ is the Riemannian
manifold ($B\times F, g),$ where
$$g = g_B+ f^2g_F.$$
\noindent More explicitly, if $U$ is tangent to $M = B\times_f F$ at
$(p, q),$ then
$$\|U\|^2 = \|d\pi_1U\|^2 + f^2(p)\|d\pi_2U\|^2$$
\noindent where $\pi_i (i=1,2)$ are the canonical projections of
$B\times F$ on $B$ and $F$, respectively.\\

\parindent=8mm They have given this important result for warped products
$$\nabla_UV = \nabla_VU = (U ln f)V\eqno(1.1)$$
\noindent for any vector fields $U$ tangent to $B$ and $V$ tangent
to $F$.\\
\par If the manifolds $M_\theta$ and $M_\perp$ are slant and
anti-invariant submanifolds respectively of a cosymplectic manifold
$\bar M$, then their warped products are
\par $(a)~M_\perp \times_f M_\theta, $
\par $(b)~M_\theta \times_f M_\perp. $\\

\parindent=8mm In the sequel, we call the warped product submanifold $(a)$ as warped product hemi-slant submanifold and the warped product $(b)$ as hemi-slant warped product submanifold.\\

\parindent=8mm Recently, K. A. Khan et. al. [5] studied warped product semi-slant submanifolds in
cosymplectic manifolds and proved that there does not exist warped
product submanifold of the type $M_1 {\times}_f M_2$ of cosymplectic
manifolds $\bar M$ where $M_1$ and $M_2$ are any Riemannian
submanifolds of $\bar M$ with $\xi$ tangential to $M_2$ other than
Riemannian product. In [9] Siraj Uddin et. al. studied warped
product submanifolds with slant factor and showed that warped
product submanifold of the type $M_1 {\times}_f M_2$ of cosymplectic
manifolds $\bar M$, such that $\xi \in TM_1$, where $M_1$ is totally
real submanifold and $M_2$ is proper slant submanifold of $\bar M$
are simply Riemannian product. So for this case we have established
an inequality for the squared norm of the second fundamental form
with constant $\phi-$sectional curvature for cosymplectic manifolds.

\parindent=8mm On the other hand, warped product submanifold of the type $M_\theta {\times}_f M_\perp$ of
cosymplectic manifolds $\bar M$, such that $\xi \in TM_\theta$,
where $M_\theta$ is proper slant submanifold and $M_\perp$ is
totally real submanifold of $\bar M$, we have established an
inequality for such type of submanifolds.\\

\section{Preliminaries}

Let $\bar M$ be an almost contact metric manifold and let
$\phi,\xi,\eta,g$ be it's almost contact metric structure. Thus
$\bar M$ is $(2n+1)$-dimensional differential manifold and $\phi$,
$\xi$, $\eta$, $g$ are respectively , a $(1,1)$-tensor field, a
vector field, a $1$-form, a Riemannian metric on $\bar M$ such that
$$\phi^2=-I+\eta\otimes\xi,~~~~ \phi\xi=0,~~~~\eta(\phi)=0,~~~~\eta(\xi)=1,~~~~ \eta(X)=g(X,\xi) \eqno (2.1)  $$
$$g(\phi X,\phi Y)=g(X,Y)-\eta(X)\eta(Y),~~~~g(\phi X,Y)=-g(X,\phi Y) \eqno(2.2)  $$
Here and in the sequel, $X,Y,Z,...$ always denote arbitrary vector
fields on $\bar M$ if it is not otherwise stated. The fundamental
$2$-form $\Phi$ of $\bar M$ is defined by $\Phi (X,Y)=g(\phi X,Y)$.

\parindent=8mm $\bar M$ is said to be almost cosymplectic if the forms $\eta$ and $\Phi $ are closed, that is,
$d\eta = 0$ and $d\Phi =0$, $d$ being the operator of the exterior
differentiation of differential forms [7]. If $\bar M$ is almost
cosymplectic and its almost contact structure $(\phi,\xi, \eta)$ is
normal, then $\bar M$ is called cosymplectic. It is well known that
a necessary and suffìcient condition for $\bar M$ to be cosymplectic
is that $\bar\nabla \phi$ vanishes identically, where $\bar \nabla$
is the Levi-Civita connection on $\bar M$. A plane section
$\sigma\subset T_x(\bar M)$ is a $\phi-section$ if $\sigma$ is
spanned by $\{u,\phi_x u\} $, for some $u\in T_x(\bar M)$. If we
restrict the $\phi-planes$ by a point function then the Riemannian
sectional curvature $(of (\bar M, g))$ is the $\phi-${\it sectional
curvature}. Now, let $\bar M(c)$ be a cosymplectic manifold of
constant $\phi-${\it sectional curvature} $c$. Then the curvature
tensor $\bar R$ of $\bar M(c)$ is given by
$$\bar R(X,Y,Z,W)=\frac{c}{4}\{g(\phi Y,\phi Z)g(X,W)-g(\phi X,\phi Z)g(Y,W) \leqno(2.3) $$
$$~~~~~~~~~~~~~~~~~~~~+\eta(Y)\eta(W)g(X,Z)-\eta(X)\eta(W)g(Y,Z)$$
$$~~~~~~~~~~~~~~~~~~~~~~~~~~~~~~~~~~~~~~~~~~~~~~~~~~~~~+g(\phi Y,Z)g(\phi X,W)-g(\phi X,Z)g(\phi Y,W)+2g(X,\phi Y)g(\phi Z,W) \} $$
\noindent for any $X,Y,Z,W \in \bar M(c) $.\\

\parindent=8mm Let $M$ be a real $m$-dimensional submanifold of $\bar M$. We shall need the Gauss-Weigarten formulae
$$\bar \nabla_X Y=\nabla_X Y+ h(X,Y),~~\bar \nabla_X V=-A_V X+ \nabla_X^\perp V, \leqno(2.4)$$
for any $X,Y\in TM$ and $V\in T^\perp M,$ where $\nabla^\perp$ is
the connection on the normal bundle $T^\perp M$, $h$ is the second
fundamental form and $A_V$ is the Weingarten map associated with the
vector field $V\in T^\perp M$ as
$$g(A_V X,Y)=g(h(X,Y),V).\leqno(2.5) $$

\parindent=8mm We denote by $\bar R$ and $R$ the curvature tensor
fields associated with $\bar \nabla$ and $\nabla$, respectively. We
recall the equation of Gauss and Codazzi
$$\bar R (X,Y,Z,W)=R(X,Y,Z,W)+g(h(X,Z),h(Y,W))-g(h(X,W),h(Y,Z)),  \leqno(2.6) $$
$$\bar R (X,Y,Z,V)=g(({\nabla^h}_X h)(Y,Z),V)-g(({\nabla^h}_Y h)(X,Z),V), \leqno(2.7)  $$
\noindent for any $X,Y,Z,W \in TM$ and $V\in T^\perp M$, where
$(\nabla^h) h$ is the covariant derivative of the second fundamental
form given by
$$ ({\nabla^h}_X h)(Y,Z)=\bar\nabla_X h(Y,Z)-g(\nabla_X Y,Z)-g(Y,\nabla_X Z), \leqno(2.8) $$
\noindent for all $X,Y,Z \in TM$. The second fundamental form $h$
satisfies the classical Codazzi equation (according to [6]) if
$$ (\nabla_X h)(Y,Z)=(\nabla_Y h)(X,Z). $$

\parindent=8mm Let $p\in M$ and $\{e_1, ..., e_m, ..., e_{2m+1}\}$ an orthonormal basis of the tangent
space $T_p \bar M(c)$, such that $e_1, ..., e_m$ are tangent to $M$
at $p$. We denote by $H$ the mean curvature vector, that is
$$H(p)=\frac{1}{m}\sum\limits_{i=1}\limits^m h(e_i,e_i). $$
\parindent=8mm Also, we set
$${h^r}_{ij}=g(h(e_i,e_j),e_r),~~~~i,j\in\{1,...,m\},r\in\{m+1,...,2m+1\}. $$
and
$$ \|h\|^2=\sum\limits_{i,j=1}\limits^m g((h(e_i,e_j),h(e_i,e_j)).$$
\parindent=8mm A submanifold $M$ is totally geodesic in $\bar M$ if $h = 0$, and minimal if $H = 0$.

\section{Hemi-slant submanifolds}
\parindent=8mm Throughout the section $M$ is a hemi-slant submanifold of
an almost contact metric manifold $\bar{M}$. Now in this section we
shall discuss hemi-slant submanifolds of cosymplectic manifolds.
More precisely, we will study integrability of the distributions of
$M$ and of the immersion of their leaves in $M$ or $\bar M$.\\

\noindent {\bf Definition 3.1. [8]} {\it  A submanifold $M$ is said
to be a hemi-slant submanifold of an almost contact metric manifold
$\bar{M}$, if there exist two orthogonal distributions $D^\perp$ and
$D_\theta$ on $M$ such that}
\begin{enumerate}
 \item[{(i)}] $T_M = D^\perp \oplus D_\theta \oplus <\xi>$.
 \item[{(ii)}] The distribution $D^\perp$ is anti-invariant $i.e.,$ $\phi D^\perp \subseteq T^{\perp}M$.
 \item[{(iii)}] The distribution $D_\theta$ is slant with slant angle $\theta \neq
\pi/2$
\end{enumerate}
\noindent from the definition it is clear that if $\theta =0$, the
hemi-slant submanifold become semi-invariant submanifold.\\
\parindent=8mm Suppose $M$ to be a hemi-slant submanifold of an almost
contact metric manifold $\bar M$. Then, for any $X \in TM,$ put
$$X= P_1 X+P_2X+\eta(X) \xi \leqno(3.1) $$
\noindent where $P_i=(i=1,2) $ are projection maps on the
distributions $D^\perp$ and $D_\theta$. Now operating $\phi$ on both
sides of equation (3.1), we have
$$\phi X=NP_1X+TP_2X+NP_2X \leqno(3.2) $$
\noindent it is easy to see that $TX=TP_2X$, $NX=NP_1X+NP_2X $,
$\phi P_1X=NP_1X$, $TP_1X=0$ and $TP_2X \in D_\theta$. Also we put
$$\phi V=BV+CV \leqno(3.3)  $$
for any $V\in T^\perp M$, where $BV$ is the tangent part of $\phi V$
and $CV$ is the normal part of $\phi V$. We define three tensor
fields $\psi_{\perp}: TM\longrightarrow T^\perp M$,
$\omega_{\theta}:TM\longrightarrow TM$ and
$\kappa_{\theta}:TM\longrightarrow T^\perp M$ by $\psi X=NP_{1}X $,
$\omega X=TP_{2} X$ and $\kappa X=NP_2 X$ respectively, for any
$X\in TM$. Now by using all the above equations and the equations of
Gauss and Weigarten for the immersion of $M$ in $\bar M$, we obtain
following lemmas which play
an important role in working out new results.\\

\noindent {\bf{Lemma 3.1.}} {\it Let $M$ be a hemi-slant submanifold
of a cosymplectic manifold $\bar M$. Then}
$$\nabla_X \omega Y-A_{\psi Y} X-A_{\kappa Y} X=\psi \nabla_X Y+\omega \nabla_X Y -Bh(X,Y) \leqno(3.4) $$
$$h(X,\omega Y)+{\nabla^\perp}_X \psi Y+{\nabla^\perp}_X \kappa Y=\kappa \nabla_X Y+Ch(X,Y) \leqno(3.5)  $$
$$\eta(\nabla_X \omega Y)=\eta(A_{\psi Y} X)+\eta(A_{\kappa Y} X) \leqno(3.6) $$
\noindent for any $X,Y \in TM$.\\

\noindent{\it Proof.} For any $X,Y \in TM$, from structure equation
we have
$$ (\bar\nabla_X \phi)Y=\bar\nabla_X \phi Y-\phi \bar\nabla_X Y=0  $$
\noindent from (3.1), we have
$$ \bar\nabla_X \phi(P_1 Y+P_2 Y+\eta(Y)\xi )-\phi\bar\nabla_X Y=0 $$
using (3.2), we get
$$\bar\nabla_X (NP_1 Y+TP_2 Y+NP_2 Y)-\phi\bar\nabla_X Y=0  $$
putting the values of tensor fields
$$\bar\nabla_X \psi Y+\bar\nabla_X \omega Y+\bar\nabla_X \kappa Y=0  $$
using Gauss and Weigarten formulae and (3.1), we get
$$-A_{\psi Y} X+{\nabla^\perp}_X \psi Y+\nabla_X \omega Y+h(X,\omega Y)-A_{\kappa Y}X+{\nabla^\perp}_X \kappa Y $$
$$~~~~~~~~~~~~-\psi \nabla_X Y-\omega \nabla_X Y+\kappa \nabla_X Y-B h(X,Y)-Ch(X,Y)=0 $$
on equating tangential and normal parts, we obtain (3.4) and (3.5).

\noindent {\bf{Lemma 3.2.}} {\it Let $M$ be a hemi-slant submanifold
of a cosymplectic manifold $\bar M$. Then}
\begin{enumerate}
\item[{(i)}]$A_{\phi Z} W=A_{\phi W}Z$ for all $W,Z \in D^\perp$,
\item[{(ii)}] $ [Z,\xi] \in D_\perp$ and $[X,\xi]\in D_\theta$ for all $Z\in D_\perp$ and $X\in D_\theta$,
\item[{(iii)}]$ g([U,V],\xi)=0 $ for all $U,V \in D^\perp\oplus D_\theta$.
\end{enumerate}

\noindent{\it Proof.} The proof is straightforward and may be
obtained
by using structure equation with equations (2.4) and (2.5) .\\

\noindent {\bf{Theorem 3.1.}} {\it Let $M$ be a hemi-slant
submanifold of a cosymplectic manifold
$\bar M,$ then the anti-invariant distribution $D_\perp$ is integrable.}\\

\noindent{\it Proof.} For any $Z,W \in D_\perp$ and $Z\in D_\theta$
by using equation (3.1)
$$g([Z,W],TP_2 X)=-g(\phi[Z,W],P_2 X)  $$
\noindent u sing Structure equation and (2.4), we get
$$g([Z,W],TP_2 X)=g(A_{\phi Z} W-A_{\phi W} Z,P_2 X)  $$
\noindent the integrability of distribution follows from Lemma
(3.2).$~~~~~~~~~~~~~~~~~~~~~~~~~~~~~~~~~~~~~~~~\square$\\

\noindent {\bf{Theorem 3.2.}} {\it Let $M$ be a hemi-slant
submanifold of a cosymplectic manifold $\bar M,$ then the slant
distribution $D_\theta$ is integrable
$$ h(X,TY)-h(X,TY)+{\nabla_X}^\perp NY-{\nabla_Y}^\perp NX$$
lies in $ND_\theta$ for each $X,Y \in D_\theta$.}\\

\noindent{\it Proof}. For any $Z\in D_\perp$, making use of
equations, we obtain
$$g(N[X,Y],NZ)=g(h(X,TY)-h(Y,TX)+\nabla _X^\perp NY-\nabla_Y^\perp NX,NZ) $$
\noindent The result follows on using the fact that $ND^\perp$ and
$ND_\theta$ are mutually perpendicular.
\parindent=8mm $~~~~~~~~~~~~~~~~~~~~~~~~~~~~~~~~~~~~~~~~~~~~~~~~~~~~~~~~~~~~~~~~~~~~~~~~~~~~~~~~~~~~~~~~~~~~~~~~~~~~~~~~~~~~~~\square$\\

\noindent {\bf{Theorem 3.3.}} {\it Let $M$ be a hemi-slant
submanifold of a cosymplectic manifold $\bar M,$ then}
\begin{enumerate}
\item[{(i)}] The leaves of the distribution $D_\perp$ are totally geodesic in $M$ if and only if $g(h(D_\perp, D_\perp),
ND_\theta)=0$.
\item[{(ii)}] The leaves of the distribution $D_\theta$ are totally geodesic in $M$ if and only if $g(h(D_\perp,
D_\theta),ND_\theta)=0$.
\end{enumerate}
\noindent{\it Proof.} (i) By assumption $g(\nabla_ZW, X)=0$  and
$g(\nabla_Z W,\xi)=0 $ for each $Z,W \in D_\perp$ and $X\in
D_\theta$, therefore
$$g(\bar\nabla_Z \phi W,\phi X)=0 $$
on using Gauss formula
$$g(h(X,\phi Y),NZ)=0 $$
Result follows from above equation.\\
(ii) Again by assumption $g(\nabla_X Y, Z)=0$  and $g(\nabla_X
Y,\xi)=0 $ for each $X,Y \in D_\theta$ and $Z\in D_\perp$, therefore
using (2.4), (2.5) and structure equation, we get
$$g(h(X,Z),NW)=0 $$
and we complete the theorem.$~~~~~~~~~~~~~~~~~~~~~~~~~~~~~~~~~~~~~~~~~~~~~~~~~~~~~~~~~~~~~~~~~~~~~~~~~~~~~~~~~~~~~~~~~~~~~~~\square$\\

\noindent {\bf{Theorem 3.4.}} {\it Let $M$ be a hemi-slant
submanifold of a cosymplectic manifold $\bar M,$ then $M$ is
hemi-slant product if and only if
$$\nabla_Z W\in D_\perp \leqno(3.7) $$
\noindent for any $Z,W \in D_\perp$}.\\
\noindent{\it Proof.} Suppose $M$ is a hemi-slant product locally
represented by $M_1 \times M_2$. Then $M_1$ and $M_2$ are totally
geodesic in $M$ then
$$\nabla_Z W={\nabla^1}_Z W \in D_\perp \leqno(3.8) $$
for any $Z,W \in D_\perp$
$$\nabla_X Y={\nabla^2}_X Y \in D_\theta \leqno(3.9) $$
for any $X,Y \in D_\theta$

\noindent where $\nabla^1$ and $\nabla^2$ are the Riemannian
connections on $M_1$ and $M_2$ respectively. Again using (3.1)
$$g(\nabla_X Z,Y)=-g(Z,\nabla_X Y)=0 \leqno(3.10)  $$
for any $X,Y \in D_\theta$ and $Z\in D_\perp$. Thus from (3.8) and
(3.10) it follows that (3.7) holds.

{\it Conversely} By using the fact that $M_1$  and $M_2$ are totally
geodesic we get the result.$~~~~~~~~~~~~~~~~~~~~~~~~~~~~~~~~~~~~~~~~~~~~~~~~~~~~~~~~~~~~~~~~~~~~~~~~~~~~~~~~~~~~~~~~~~~~~~~\square$\\

\noindent {\bf{Theorem 3.5.}} {\it Let $M$ be a hemi-slant
submanifold of a cosymplectic manifold $\bar M$. Then $M$ is
hemi-slant product if and only if its second fundamental form
satisfies
$$Bh(X,Z)=0 \leqno(3.11) $$
$$h(X,\phi Y)=Ch(X,Y),  \leqno(3.12)  $$
\noindent for any $X\in TM$ and $Y\in D_\perp$.}\\

\noindent{\it Proof}. From (3.4) and (3.5) it follows that
$$\nabla_X \phi Z=\psi(\nabla_X Z)+Bh(X,Z) \leqno(3.13) $$
and
$$h(X,\phi Z)=Ch(X,Z)+\kappa (\nabla_X Z), \leqno(3.14)  $$
\noindent for any $X\in TM$ and $Z\in D_\perp$. Thus our assertion
follows from (3.13) and (3.14) by means of Theorem (3.4).\\

\parindent=8mm Now, using the formulas of Gauss and Weigarten, we
obtain
$$g(A_{\phi X} Z,Y)=-g(Bh(Z,Y),X), \leqno(3.15)  $$
\noindent for any $Z\in D_\perp$, $Y\in TM$ and $X\in D_\theta$. $~~~~~~~~~~~~~~~~~~~~~~~~~~~~~~~~~~~~~~~~~~~~~~~~~~~~~~~~~~~~~~~~~~~\square$\\

\parindent=8mm Then, by using above equation and Theorem (3.5), we obtain the following corollary\\

\noindent {\bf{Corollary 3.1.}} {\it Let $M$ be a hemi-slant
submanifold of a cosymplectic manifold $\bar M$. Then the following
assertions are equivalent to each other:}
\begin{enumerate}
\item[{(i)}]$M$ is hemi-slant product;
\item[{(ii)}] the fundamental tensors of Weingarten satisfy $A_{\phi X}Z=0,$ for any $X\in D_\theta$ and $Z\in
D_\perp$;
\item[{(iii)}] the second fundamental form of $M$ satisfies
$h(Y,\phi Z)=\phi h(Y,Z)$, for any $Z\in D_\perp$ and $Y\in TM$.
\end{enumerate}

\parindent=8mm To close this section, we recall that bisectional
curvature of a cosymplectic manifold $\bar M$ is defined by
$$ S(X,Y)=\bar R(X,\phi X, \phi Y,Y), \leqno (3.16)  $$
\noindent where $X$ and $Y$ are unit vector fields.\\

\section{Some Inequalities for hemi-slant products}

In this section, we obtain an equality for the squared norm of the
second fundamental form in terms of constant $\phi-${\it sectional
curvature} for hemi-slant products in cosymplectic manifolds $\bar
M$. Also, we have proved an inequality for hemi-slant warped product
submanifolds in a cosymplectc manifold $\bar M$ and considered the
equality case.\\

\noindent {\bf{Theorem 4.1.}} {\it Let $M$ be a hemi-slant product
of a cosymplectic manifold $\bar M$. Then:
$$\frac{1}{2} S(Y,Z)=\|h(Y,Z)\|^2-1 \leqno (4.1) $$
for any unit vector fields $Y\in D_\perp$, $Z\in D_\theta$.}\\

\noindent{\it Proof.} By using (2.4), (2.7), Corollary 3.1 and
structure equation, we get
$$\bar R(Y,\phi Y,Z,\phi Z)=g(({\nabla^h}_Y h)(\phi Y,Z)-g(({\nabla^h}_{\phi Y} h)(Y,Z),\phi Z)  $$
\noindent for $Y\in D_\perp$ and $Z\in D_\theta$. From (2.8), we
have
$$\bar R(Y,\phi Y,Z,\phi Z)=g(\bar\nabla_Y h(\phi Y,Z)-h(\nabla_Y \phi Y,Z)-h(\phi Y,\nabla_Y Z),\phi Z)$$
$$~~~~~~~~~~~~~~~~~~~~~~~-g(\bar\nabla_{\phi Y} h(Y,Z)-h(\nabla_{\phi
Y}Y,Z)-h(Y,\nabla_{\phi Y} Z),\phi Z)  $$ \noindent as $g$ is
Riemannian metric, we arrive at
$$\bar R(Y,\phi Y,Z,\phi Z)=2-g(h(\phi Y,Z),\phi \nabla_Y Z)-g(h(\phi Y,Z),\phi h(Y,Z)) $$
$$~~~~~~~~~~~~~~~~~~~~~~~~+g(h(Y,Z),\phi\nabla_{\phi Y} Z)+g(h(Y,Z),\phi h(\phi Y,Z)).  $$
\parindent=8mm Now, using (3.12) and (3.16), we get
$$ \frac{1}{2}S(Y,Z)=\|h(Y,Z)\|^2-1 $$
\noindent we complete the proof.$~~~~~~~~~~~~~~~~~~~~~~~~~~~~~~~~~~~~~~~~~~~~~~~~~~~~~~~~~~~~~~~~~~~~~~~~~~~~~~~~~~~~~~~~~~~~~~~\square$\\

\parindent=8mm Now, we will prove a condition for existence of hemi-slant products in cosymplectic manifolds in terms of $\phi-${\it sectional
curvature} $c$.\\

\noindent {\bf{Theorem 4.2.}} {\it There exist no proper hemi-slant
products in a cosymplectic manifold $\bar M(c)$ with $c\geq 2$.}\\

\noindent{\it Proof.} From (2.3) and (4.1), we get
$$\|h(X,Z)\|^2=\frac{c}{2} \|X\|^2 \|Z\|^2 \leqno(4.2) $$
\noindent for any unit vector fields $X\in D_\theta$ and $Z\in
D_\perp$. We get $c\geq2$ and this completes the proof.$~~~~~~~~~~~~~~~~~~~~~~~~~~~~~~~~~~~~~~~~~~~~~~~~~~~~~~~~~~~~~~~~~~~~~~~~~~~~~~~~~~~~~~~~~\square$\\

\parindent=8mm Let $\{\xi=e_0,e_1,....,e_p,E_1,E_2,...,E_n,\phi e_1,\phi e_2,...,\phi e_p,\phi E_1,\phi E_2,...,\phi E_n
\}$ be an orthonormal basis of $\bar M$ such that,
$\{\xi=e_0,e_1,....,e_p,E_1,E_2,...,E_n\}$ are tangent to M. Such
that $\{\xi=e_0,e_1,....,e_p\}$ form an orthonormal frame of
$D_\perp$ and $\{E_1,E_2,...,E_n\}$ form an orthonormal frame of
$D_\theta$ (where $n$ is even). We can take $\{\phi e_1,\phi
e_2,...,\phi e_p,\phi E_1,\phi E_2,...,\phi E_n \}$ as orthonormal
frame of $T^\perp \bar M$.\\

\noindent {\bf{Theorem 4.3.}} {\it Let $\bar M$ be a proper
hemi-slant products in a cosymplectic manifold $\bar M(c)$, with
$c\geq 2$. Then}
$$\|h\|^2\geq np(c-2). \leqno(4.3) $$

\noindent{\it Proof.} By adopting above frame, we get
$$\|h\|^2=\sum\limits_{i,j=1}\limits^p \|h(e_i,e_j)\|^2+\sum\limits_{\alpha,\beta=1}\limits^n \|h(E_\alpha,E_\beta)\|^2+2\sum\limits_{i,\alpha}\|h(e_i,E_\alpha)\|^2 $$
$$\|h\|^2 \geq 2\sum\limits_{i,\alpha}\|h(e_i,E_\alpha)\|^2 $$
$$\|h\|^2 \geq np(c-2) $$
and so, we have (4.3).$~~~~~~~~~~~~~~~~~~~~~~~~~~~~~~~~~~~~~~~~~~~~~~~~~~~~~~~~~~~~~~~~~~~~~~~~~~~~~~~~~\square$\\

\parindent=8mm We prove following lemma to get the inequality for warped product submanifolds.\\

\noindent {\bf{Lemma 4.1.}} {\it Let $M_\theta \times_f M_\perp$,
such that $\xi\in TM_\theta$ be a hemi-slant product submanifold of
a cosymplectic manifold $\bar M$, then we have}
$$g(h(X,Z),FZ)=X (ln f)g(Z,FZ) \leqno(4.4)  $$
$$g(h(X,Y),FZ)=0 \leqno(4.5)  $$
\noindent {\it for $X,Y\in D_\theta$ and $Z,W \in D_\perp$ }.\\

\noindent {\it Proof.} From Gauss formula $g(h(X,Z),FZ)=g(\bar\nabla_Z X,FZ)=g(\nabla_Z X, FZ)=X (ln f)g(Z,FZ) $.\\
\par Similarly again by Gauss formula
$$g(h(X,Y),\phi Z)=0. $$
This proves the Lemma
completely.$~~~~~~~~~~~~~~~~~~~~~~~~~~~~~~~~~~~~~~~~~~~~~~~~~~~~~~~~~~~~~~~~~~~~~~~~~~~~~~~~~\square$\\

\parindent=8mm Now, using the above theorem we have the following main result.
We are going to obtain an inequality for the squared norm of the
second fundamental form in terms of the warping function for
hemi-slant warped product submanifold $M$ of a cosymplectic manifold $\bar M$.\\

\noindent {\bf{Theorem 4.4.}} {\it Let $ M=M_\theta \times_f M_\perp
$ be a hemi-slant warped product submanifold of cosymplectic
manifold $\bar M$. Then, the norm of the second fundamental form of
$M$ satisfies
$$\|h\|^2\geq 2c\|\nabla ln f\|^2, $$
where $\nabla (ln f)$ is the gradient of $ln f$ and $c $ is the
dimension of $M_\perp$. If the equality sign holds then $M_\theta$
and $M_\perp$ are totally geodesic submanifolds of $M$}.\\

\noindent{\it Proof.} By adopting the frame
$$ \|h\|^2=\|h(D_\theta,D_\theta)\|^2+2\|h(D_\theta, D_\perp)\|^2+\|h(D_\perp, D_\perp)\|^2,  $$
\noindent using (4.4) and (4.5), we get
$$ \|h\|^2=\|h(D_\theta,D_\theta)\|^2+2\sum\limits_{\alpha=1,i=1}\limits^{n,p}\|h(E_\alpha, e_i)\|^2+\|h(D_\perp, D_\perp)\|^2,  \eqno(4.6) $$
The first part of theorem follows from above inequality. If the
equality sign holds, then from equation (4.6), we get
$h(D_\theta,D_\theta) = 0, h(D_\perp,D_\perp) = 0$. Since $M_\theta$
and $M_\perp$ are totally geodesic submanifolds of $M$. This
complete the proof.
$~~~~~~~~~~~~~~~~~~~~~~~~~~~~~~~~~~~~~~~~~~~~~~~~~~~~~~~~~~~~~~~~~~~~~~~~~~~~~~~~~~~~~~~~~~~~~~~~~~~~~~~~~~~~~~~~~~\square$\\

\vspace{.2cm}
\noindent{\it Author's addresses:}\\

\vspace{-.2cm}
\noindent{\small {Khushwant Singh and S. S. Bhatia}\\
School of Mathematics and Computer Applications\\
Thapar University, Patiala-147 004, INDIA\\
\noindent {\it E-mail}: {\tt khushwantchahil@gmail.com, ssbhatia@thapar.edu}}\\

\end{document}